\documentclass{amsart}

\newtheorem{theorem}{Theorem}[section]
\newtheorem{lemma}[theorem]{Lemma}

\usepackage[table,xcdraw]{xcolor}

\theoremstyle{definition}
\newtheorem{definition}[theorem]{Definition}

\theoremstyle{remark}
\newtheorem{remark}[theorem]{Remark}

\numberwithin{equation}{section}
\usepackage{url}
\usepackage{amssymb}

\usepackage{listings}
\usepackage{xcolor} 

\usepackage{hyperref}
\begin{document}
	
	\title{Sharing Nim and Enumeration of Nim Characteristics}
	
	\author{Donghyun (Ethan) Kim}
	\address{Gyeonggi Suwon International School, 451 YeongTong-Ro, YeongTong-Gu, Suwon-Si, Gyeonggi-Do, Republic of Korea}
	
	\email{ethank11k@gmail.com}

	

	\dedicatory{{\normalfont Gyeonggi Suwon International School, 451 YeongTong-Ro, YeongTong-Gu, Suwon-Si, Gyeonggi-Do, Republic of Korea}}

	\begin{abstract}
		In this paper, we introduce and examine a variant of the game of Nim (Sharing Nim), where players can either remove or transfer objects from 1 pile to another. The only restriction is that players may not transfer objects from a pile of greater size to a pile of smaller size. We also find new methods of enumerating characteristics of Nim and Sharing Nim, including the number of zero nim positions.
	\end{abstract}
	
	\maketitle
	
	
	
	
	\section{Introduction}

	Impartial games have been studied for decades as a part of the field of combinatorial game theory. The Game of Nim, which originates from ancient China, was studied by Charles L. Bouton, who developed the theory behind the game in 1901. Before examining the game of Nim and a variant of Nim, we will introduce some standard definitions and notations relevant to the study of impartial games.

\section{Some Definitions and Notation}

\begin{definition}
	An impartial game is a 2 player game in which the moves available to each player depends on the position of the game and not on which player's turn it is. 
\end{definition}

There are two versions of impartial games.

\begin{definition}
	In an impartial game, Normal Play refers to a version of the game in which the player unable to move loses.
\end{definition}

\begin{definition}
	Misère play refers to a version of the game in which the player unable to move wins.
\end{definition}

There are two impartial game positions as well.

\begin{definition}
	A P-position is defined as a game position in which the Previous player (the player who just moved) secures a win.
\end{definition}

\begin{definition}
	A N-position is defined as a game position in which the Next player secures a win.
\end{definition}

\section{Game of Nim}

In this section, we will examine the game of Nim. Note that much of this section originates from Bouton (\cite{A}) 

\begin{definition}
	In the game of Nim, 2 players take turns removing any number of objects from a single pile (there is no limit on the total number of piles). The player who takes the last object is the winner of the game.
\end{definition}

\begin{definition}
	The Nim-sum $\bigoplus$ is the Exclusive-Or (XOR) sum of all heaps. The bit-wise exclusive-or operation is given by writing the numbers in binary and carrying out addition without carry.
\end{definition}

\begin{definition}
	The bit-wise AND ($\&$) operator between 2 integers in binary form compares each bit of the 2 integers. If both bits are 1, the resulting bit is set to 1. In all other cases, the resulting bit is set to 0. 
	\label{sec:lemma-3.3}
\end{definition}

\begin{lemma}
	The sum of 2 positive integers $a$ and $b$ can be given by the following.$\\$
	\label{sec:lemma-3.4}
	$a + b = (a \bigoplus b) + (a \& b) \times 2$
\end{lemma}

Proof. Considering the sum bit-by-bit, the bit-wise XOR operator gives the corresponding bit for all cases of addition. $\\$
$0 + 0 = 0 \\ 0 + 1 = 1 \\ 1 + 0 = 1 \\ 1 + 1 = 0 \\$

In the case that carrying of bits occur (4th case from above), the bit-wise AND operator gives 1 for the resulting bit. To carry the 1, the bit can be shifted to the next position on the left (multiply by 2). 

\begin{lemma} {\normalfont(\cite{A})}
	If after a player's turn the Nim-sum equals 0, the next player must change it.
	\label{sec:lemma-3.5}
\end{lemma}

Proof. Let the initial number of objects in the piles be $x_1,x_2,...,x_n$ and the nim-sum of the piles be $S = x_1 \bigoplus x_2 \bigoplus ... \bigoplus x_n$. Let the number of objects in the piles after the next player goes be $y_1,y_2,...,y_n$ and nim-sum of the piles be $K = y_1 \bigoplus y_2 \bigoplus ... \bigoplus y_n$. $\\$

Since $S=0$ and only 1 pile is changed, after the next player goes, some $x_m \neq y_m$ and the other piles $x_n = y_n$ for $m \neq n$. $\\$

With this in mind, an equation can be made. 

$K = 0 \bigoplus K \\ = S \bigoplus S \bigoplus K \\ = S \bigoplus (x_1 \bigoplus x_2 \bigoplus ... \bigoplus x_n) \bigoplus (y_1 \bigoplus y_2 \bigoplus ... \bigoplus y_n) \\ = S \bigoplus (x_1 \bigoplus y_1) \bigoplus (x_2 \bigoplus y_2) \bigoplus ... \bigoplus (x_m \bigoplus y_m) \\ = S \bigoplus x_m \bigoplus y_m$

As $S=0$, the value of $K$ will always be non-zero (the value of $x_m \bigoplus y_m$ is non-zero).

\begin{lemma} {\normalfont(\cite{A})}
	If after a player's turn the Nim-sum is not equal to 0, the next player can always make it become 0.
	\label{sec:lemma-3.6}
\end{lemma}

Proof. The proof for Lemma 3.4 follows directly from the sum $K$ in Lemma 3.3. Let $d$ be the position of the most significant bit in $S$ and choose a pile $x_m$ such that its most significant bit is in position $d$. Then the value of pile $y_m$ can be altered to $y_m = S \bigoplus x_m$ by removing $x_m - y_m$ objects from the pile. The new Nim-sum becomes:$\\$

$K = S \bigoplus x_m \bigoplus y_m \\ = S \bigoplus x_m \bigoplus x_m \bigoplus S \\ = S \bigoplus S \bigoplus x_m \bigoplus x_m \\ = 0$ 

\begin{theorem} {\normalfont(\cite{A})} 
	If the starting game position has a Nim-sum not equivalent to 0, the first player to go will win the game. If the starting game position has a Nim-sum  equivalent to 0, the first player to go will lose the game. 
	\label{sec:lemma-3.7}
\end{theorem}

Proof. If the starting game position has a Nim-sum not equivalent to 0, the first player can make the Nim-sum equal to 0 (by Lemma 3.4). The second player will then disturb the 0 Nim-sum (by Lemma 3.3) and the first player will return it to 0. This process will repeat until the first player reaches a 0 Nim-sum with no more objects remaining. $\\$

However, if the starting game position has a Nim-sum equivalent to 0, the second player will have the winning strategy (following the same process as above).

\section{Sprague-Grundy Theorem}

In this section, we will take a look at the Sprague-Grundy Theorem. Before we move on to the theorem, we will place games on graphs as shown below. 

A game can be represented on a graph $G = (X, F)$, in which:

\begin{itemize}
	\item $X$ is the set of all game positions
	\item For each $x \in X$, the function $F$ gives a choice (subset) of possible $x$'s (positions) to move to, called followers.
	\item The starting game position is $x_0 \in X$
	\item From position $x$, a player chooses $y \in F(X)$
	\item Once $F(x)$ is empty, the next player to go loses (terminal position).
	\item The graph is progressively bounded. That is, the graph is finite and does not contain cycles.
\end{itemize}

\begin{definition} 
	The minimum excludent function (mex) is a function that gives the smallest non-negative integer that is not found among the elements of a set $N$.
	
	mex(N) = min\{$x \in \mathbb{Z}_{\geq 0} : x \notin N$\}
	
\end{definition}

\begin{definition} (\cite{B},\cite{C})
	On a game $G = (X,F)$, the Sprague-Grundy function is defined as $g: X -> \mathbb{Z}_{\geq 0}$, in which: $\\$
	
	$g(x) =$ min$ \{ j \geq 0 : j \neq g(y)$ for $y \in F(x) \}$ = mex $\{ g(y): y \in F(x) \} \\$ 
	
	This function returns the smallest non-negative integer not found in the set of Sprague-Grundy values of $x$'s followers.
	
\end{definition}

\begin{definition} (\cite{B}, \cite{C}, \cite{D})
	The graphs of impartial games may be added (disjunctive sum of 2 or more games). The disjunctive sum of any $n$ games $G_i$ is given as follows. $\\$
	
	When games $G_1 = (X_1, F_1), G_2 = (X_2, F_2),...G_n = (X_n, F_n)$ are added, $G(X,F) = G_1 + G_2 + ... + G_n$ where $X = X_1 \times X_2 \times ... \times X_n$ and $F = F_1 \times F_2 \times ... \times F_n$.
	
\end{definition}

\begin{theorem} {\normalfont(\cite{B}, \cite{C}, \cite{D})}
	The Sprague-Grundy function for a disjunctive sum of games is given by the Nim sum of the Sprague-Grundy functions of its component games. $\\$
	
	If $g_i$ is the Sprague-Grundy function of $G_i$, $i = 1,...,n$, then the Sprague-Grundy function of $G = G_1 + ... + G_n$ is given by $g(x_1,...,x_n)  = g_1(x_1) \bigoplus ... \bigoplus g_n(x_n)$.
	\label{sec:theorem-4.4}
\end{theorem}

Proof. Let $x = (x_1,...,x_n)$ be an arbitrary point of $X$ and let $b = g_1(x_1)\bigoplus...\bigoplus g_n(x_n)$. We will proceed in 2 parts. $\\$

1. For every non-negative integer $a < b$, there is a follower of $x$ that has a Sprague-Grundy value $a$. $\\$

Let $d = a \bigoplus b$ and let $k$ be the number of digits in the binary expansion of $d$. Since $a < b$, $b$ has a 1 in its $k$th digit. Also, since $b = g_1(x_1)\bigoplus...\bigoplus g_n(x_n)$, at least one $x_i$ has a $g_i(x_i)$ value with a 1 in its $k$th digit (in its binary expansion). $\\$

Without loss of generality, we can set $i=1$. Then we see that $d \bigoplus g_1(x_1) < g_1(x_1)$ and thus a move from $x_1$ to some $x_i^\prime$ with $g_1(x_1^\prime) = d \bigoplus g_1(x_1)$ exists. $\\$

Hence the move from $(x_1,x_2,...,x_n)$ to $(x_1^\prime, x_2, ..., x_n)$ is a valid move in $G$ and as a result, $g_1(x_1^\prime) \bigoplus g_2(x_2) \bigoplus ... \bigoplus g_n(x_n) = d \bigoplus g_1(x_1) \bigoplus g_2(x_2) \bigoplus ... \bigoplus g_n(x_n) = d \bigoplus b = a \\$

2. No follower of $x$ has a Sprague-Grundy value $b$. $\\$

Assume for the sake of contradiction that $(x_1,...,x_n)$ has a follower with the same Sprague-Grundy value and without loss of generality let this involve a turn in the first game. In other words, assume that $(x_1^\prime, x_2,...,x_n)$ is a follower of $(x_1,x_2,...x_n)$ and that $g_1(x_1) \bigoplus g_2(x_2) \bigoplus ... \bigoplus g_n(x_n)$. Through the cancellation property of the XOR sum, $g_1(x_1^\prime) = g_1(x_1)$. However, this is a contradiction as no position can have the same Sprague-Grundy value as a follower.

\begin{remark}
	Theorem \ref{sec:theorem-4.4} can be applied to our discussion of the game of Nim above. If each pile is thought of as an individual game of Nim where the Sprague-Grundy value is the size of the pile, then the Sprague-Grundy value of the game is the nim-sum of individual piles of objects. 
\end{remark}

	\section{Sharing Nim}
		
	In this section, we introduce a new variant of the Game of Nim: Sharing Nim. $\\$

	\begin{definition}
		\emph{Sharing Nim} is a game played with the following rules. The game starts with 3 piles of objects that each contain 1 or more objects. 2 players take turns either removing a certain number of objects from 1 pile or transferring a certain number of objects from 1 pile to another. The only restriction is that players may not transfer objects from a pile of greater size to a pile of smaller size. The game ends when no more objects remain and the last player to remove an object wins. 
	\end{definition}

	We begin by analyzing cases of the game with less than 3 piles of objects. $\\$
	
	$\textbf{Game with 1 pile:}$ $\\$
	
	If the game begins with 1 pile of objects, Player 1 can secure the win by removing all objects at once (starting position is a $N$-position). $\\$
	
	$\textbf{Game with 2 piles:}$ $\\$
	
	If the game begins with 2 piles of objects, the game can be analyzed similar to the game of Nim.

	\begin{lemma}
		The properties of the game of Nim are preserved in a game of sharing Nim with 2 piles.
		\label{sec:lemma-5.2}
	\end{lemma}

	Proof. We will proceed with the proof by showing that the transferring property of Sharing Nim preserves the properties mentioned in lemmas \ref{sec:lemma-3.5} and \ref{sec:lemma-3.6}.
	
	Because there are only 2 piles, the nim-sum is 0 only if the 2 piles are equal in size.
	
	Let both piles have $a$ objects each (nim-sum = 0). If $b$ objects ($b > 0$) are transferred between the 2 piles, the piles will have $a-b$ and $a+b$ objects each. Thus, from a position where the nim-sum equals 0, all moves will result in a nim-sum unequal to 0.

	From a position where the nim-sum is not equal to 0, the nim-sum cannot be made to equal 0 by transferring objects (objects cannot be moved from piles with more objects to less objects). As such, players will instead remove objects from the pile of greater size (to make nim-sum = 0).
	
	\begin{theorem}
		In a game of Sharing Nim with 2 piles, the winning strategy and the winner is equivalent to a game of Nim. 
	\end{theorem}

	Proof. As shown in Lemma \ref{sec:lemma-5.2}, if after a player's turn the Nim-sum equals 0, the next player must change it. Also, if after a player's turn the Nim-sum is not equal to 0, the next player can always make it become 0. Thus, following the same reasoning as Theorem \ref{sec:lemma-3.7}, if the starting game position has a Nim-sum = 0, the first player to go will lose the game. If the starting game position has a Nim-sum $\neq 0$, the first player to go will win the game. $\\$
	 
	In the next section, we will note some characteristics of 3 pile Sharing Nim. 
	
	\section{Characteristics of Sharing Nim}

	\begin{definition}
		Game positions where the nim-sum is 0 will be referred to as \emph{zero-nim positions}. Note that this definition applies for both the game of Nim and its variant, Sharing Nim.
	\end{definition}

	\begin{definition}
		From certain zero-nim positions, the 2 players can repeatedly move to positions where the nim-sum is 0. These sets of positions will be referred to as a \emph{transferable sequence}. 
		
		An example of such a sequence of positions is $(7, 11, 12), (6, 11, 13), (5, 11, 14), (4, 11, 15)$. Players can move through these positions by transferring 1 object from pile 1 to pile 3. 
	\end{definition}
	
	\begin{lemma}
		Game positions of the form "a a b" are N-positions.
	\end{lemma}

	Proof. The next player can remove all the objects from the third pile. Then the game is equivalent to a game of Nim with 2 piles (as shown in Lemma \ref{sec:lemma-5.2}); thus, the player who removed all objects from the third pile  will win the game.

	\begin{lemma}
		Certain game positions with nim-sum equivalent to 0 are of the form "a, b, a+b", where $a < b$. These positions are N-positions.
		\label{sec:lemma-6.4}
	\end{lemma}

	Proof. The next player can transfer all objects in the first pile to the second pile. Then the game is equivalent to a game of Nim with 2 piles and the player who made this move will win the game.

	\begin{lemma}
		The move where a player transfers objects from one of the 3 piles to another one of the 3 piles when the nim-sum is not zero results in a position where the nim-sum is not zero.
	\end{lemma}
	
	Proof. We consider the binary expansion of the piles. $\\$
	
	$n$th digit: 1 1 1 (subtract $2^{n-1}$ from second pile and add to third pile) $\\$
	$n$th digit after move: 1 0 0 $\\$
	
	$n$th digit: 1 0 0 (subtract $2^{n-1}$ from second pile and add to third pile) $\\$
	$n$th digit after move: 1 1 1 $\\$ 
	
	$n$th digit: 0 1 0 (subtract $2^{n-1}$ from second pile and add to third pile) $\\$
	$n$th digit after move: 0 0 1 $\\$
	
	$n$th digit: 0 0 1 (subtract $2^{n-1}$ from second pile and add to third pile) $\\$
	$n$th digit after move: 0 1 0 $\\$
	
	Clearly, the move where a player transfers objects from one of the 3 piles to another one of the 3 piles when the nim-sum is not zero results in a position where the nim-sum is again unequal to zero.

	\section{Enumeration of Such Characteristics}	
	
	In this section, we will investigate methods of counting some of the characteristics from section 6. We begin with a function that is referred to throughout this section.
	
	\begin{theorem}
		The number of 1s in a binary number $n$ is given by the function $g(n)$.
		
		$g(n) = \begin{cases}
		0 & n = 0 \\
		1 & n = 1 \\
		n \bmod 2 + g(\lfloor \frac{n}{2} \rfloor) & \text{else}
		\end{cases} \\$
		
		The function can be re-written in the following form as well.$\\$
		
		$g(n) = \sum_{i = 0}^n \left( \lfloor\frac{n}{2^i}\rfloor \right) \bmod 2.\\$

	\end{theorem}

	The next theorem is also important to the section as a whole.

	\begin{theorem}
		The number of pairs of integers $a$ and $b$ that have sum $S$ and XOR sum $X$ is given by the following function $f(S,X)$.

		$f(S,X)=\begin{cases} 0 & (S - X) \& 1 = 1 \\ 0 & \text{ith bit of } \frac{S-X}{2} \text{ and } X \text{ both equal 1} \\ 2^{g(X)-1} - 1 & S = X \\ 2^{g(X)-1} & (S-X)\&1 = 0 \end{cases}$
		\label{sec:lemma-7.2}
	\end{theorem}

	Proof. 
	Let $a + b = S$ and	$a \oplus b = X.$
	
	From Lemma \ref{sec:lemma-3.4}, $a + b = (a \oplus b) + (a \& b)\times 2$. Thus, $S = X + (a \& b)\times 2$ and $(a \& b) = \frac{(S-X)}{2}.$
	
	Now, we consider the equation bit-by-bit. There are two cases that must be considered.$\\$
	
	(1). $a_i \& b_i = 1$
	
	In this case, $a_i \oplus b_i$ must equal 0 and $a_i = b_i$, giving 1 case. If $a_i \oplus b_i = 1$, then there are no such pairs of integers $a$ and $b$.$\\$
	
	(2). $a_i \& b_i = 0$
	
	In this case, there are two different cases for the value of $a_i \oplus b_i$. 
	
	If $a_i \oplus b_i = 1$, then $a_i$ and $b_i$ can take the values (0, 1) and (1,0), giving 2 cases.
	
	If $a_i \oplus b_i = 0$, then $a_i$ and $b_i$ are both equal to 0, giving 1 case.$\\$
	
	There are also two exceptions that need to be considered. $\\$
	
	(1) If $S = X$, then 2 must be subtracted from the total number of pairs of integers. This is because two cases where $a$ or $b$ equals 0 is included in the total number and both integers must be positive. $\\$
	
	(2) If $\frac{S - X}{2}$ is not an integer, there are no solutions.$\\$
	
	From the above, we can build formulas for the number of pairs of integers. $\\$
	
	From case (2), for every 1 in the binary expansion of $X$, the number of total pairs increases by a factor of 2, giving the formula $2^k$. 
	
	From exception (1), if $S = X$, 2 needs to be subtracted from the number of total pairs, giving the formula $2^k - 2$.
	
	From case (1), if $a_i \oplus b_i = 1$ and $a_i \& b_i = 1$, then there are 0 pairs of integers. In other words, if the ith bit of $\frac{S-X}{2}$ and $X$ both equal 1, then there are 0 pairs of integers.
	
	From exception (2), if $(S - X)$ is not a multiple of 2, there are 0 pairs of integers.
	
	Finally, the total number of pairs of integers must be divided by 2 to account for over counting from case (2): pairs of integers (a,b) and (b,a) are considered distinct.

	\begin{lemma}
		The number of 3 pile zero-nim positions given a fixed pile of $X$ objects and sum of 2 other piles $S$ is given by the following function $f(S,X)$.

		$f(S,X)=\begin{cases} 0 & (S - X) \& 1 = 1 \\ 0 & \text{ith bit of } \frac{S-X}{2} \text{ and } X \text{ both equal 1} \\ 2^{g(X)-1} - 1 & S = X \\ 2^{g(X)-1} & (S-X)\&1 = 0 \end{cases}$
		\label{sec:lemma-7.3}
	\end{lemma}

	Proof. The 3 piles need to have an XOR sum of 0. Hence, the 2 piles other than the fixed pile must have a XOR sum of $X$. The problem is then equivalent to Theorem \ref{sec:lemma-7.2} of finding the number of pairs of integers that have a sum $S$ and XOR sum $X$.

	\begin{lemma}
		The number of 3 pile zero-nim positions given a sum of 3 numbers $S$ is given by the following function $q(S)$, where $a_i$ and $b_i$ are the ith bits of the sizes of the 2 other piles $a$ and $b$ given a fixed pile.
		
		$q(S) = \frac{\sum_{n=1}^{S-2} f(n,S-n)}{3}$, where $f(n,S-n)=\begin{cases} 0 & (2n - S) \& 1 = 1 \\ 0 & \text{ith bit of } \frac{2n - S}{2} \text{ and } S-n \text{ both equal 1} \\ 2^{g(S-n)-1} - 1 & S = 2n \\ 2^{g(S-n)-1} & (2n-S)\&1 = 0 \end{cases}$
		\label{sec:lemma-7.4}
	\end{lemma}

	Proof. This formula is a direct extension of the formula established in lemma \ref{sec:lemma-7.3}. To compute the number of zero-nim positions, we iterate through different values for the fixed pile and use the same function as above. $\\$
	Then, to account for over-counting, we divide the entire summation by 3: given a position $a, b, c$, all 3 numbers can take the position of the fixed pile (positions with repeating numbers do not appear as the 3 numbers need a XOR value of 0).

	\begin{lemma}
		The number of 3 pile zero-nim positions with the total number of objects less than or equal to $K$ is given by the following function $p(K)$, where $a_i$ and $b_i$ are the ith bits of the sizes of the 2 other piles $a$ and $b$ given a fixed pile.
		
		$p(K) = \frac{ \sum_{S=3}^{K} \sum_{n=1}^{S-2} f(n,S-n)}{3}$, where $f(n,S-n)=\begin{cases} 0 & (2n - S) \& 1 = 1 \\ 0 & \text{ith bit of } \frac{2n - S}{2} \text{ and } S-n \text{ both equal 1} \\ 2^{g(S-n)-1} - 1 & S = 2n \\ 2^{g(S-n)-1} & (2n-S)\&1 = 0 \end{cases}$ $\\$
		
		or, with the inclusion of $q(S)$, $p(K) = \frac{ \sum_{S=3}^{K} q(s)}{3}$
		\label{sec:lemma-7.5}
	\end{lemma}

	Proof. This formula is a direct extension of Lemma \ref{sec:lemma-7.4}. To compute the number of zero-nim positions with the total numbers of objects less than or equal to $Q$, we iterate through different values for the total number and apply the same function from above.

	\begin{lemma}
		The binary digits in positions of the form $"a, b, a + b"$ (from Lemma 6.4) are restricted to $(0, 1, 1)$, $(1, 0, 1)$, and $(0, 0, 0)$ for each digit place; ie. no carrying of digits occur in $a+b$.
		\label{sec:lemma-7.6}
	\end{lemma}

	Proof. First, we note that these 3 sets of digits ensure that the positions have a nim-sum of 0 as there are an even number of 1s in each digit.$\\$
	
	Digits $1, 1, 0$ cannot appear in positions of the form $"a, b, a + b"$. If carrying of digits occur when $a$ and $b$ are added together, then the number of 1s will be odd in the next units place. There are 3 cases to consider for the digits after $1, 1, 0$. $\\$
	
	Case 1. $(0, 1, 1) -> (0, 1, 0) \\$
	
	Case 2. $(0, 0, 0) -> (0, 0, 1) \\$
	
	Case 3. $(1, 1, 0) -> (1, 1, 1) \\$
	
	As the number of 1s in the next digit place is always odd, digits $1, 1, 0$ cannot appear in positions of the form $a, b, a + b$.
	
	\begin{theorem}
		The number of unordered positions of the form $"a, b, a + b"$ (from Lemma \ref{sec:lemma-6.4}) with two piles with $k$ digits (in binary form) is given by the formula $3^{k-1} - 2^{k-1}$.
		\label{sec:lemma-7.7}
	\end{theorem}

	Proof. The $k$th digits are restricted to $1, 0, 1$ and $0, 1, 1$ in order. However, as these 2 cases are essentially the same, we arbitrarily select the first set of digits ($1,0,1$). The remaining $k-1$ digit places can be any set of digits from the following (once again, in order): $(0, 1, 1)$, $(1, 0, 1)$, and $(0, 0, 0)$ (from Lemma \ref{sec:lemma-7.6}). $\\$
	
	Then, the total number of positions is $3^{k-1}$. However, there are $2^{k-1}$ cases where one pile has 0 objects (by filling up the remaining $k-1$ digits with only the 2nd and 3rd set). Therefore, the formula is $3^{k-1} - 2^{k-1}$.
	
	\begin{lemma}
		All position of the form $"a, b, a + b"$ will have 2 piles with the same number of digits in binary form and 1 pile with less digits.
	\end{lemma}

	Proof. From Lemma \ref{sec:lemma-7.6}, the binary digits in such positions are restricted to $(0, 1, 1)$, $(1, 0, 1)$, and $(0, 0, 0)$ in order (for each digit place). Thus, the set of digits farthest to the left must be either $(0, 1, 1)$ or $(1, 0, 1)$ (once again, in order). Hence, 2 piles will have the same number of digits and 1 pile will have less digits in binary form.
	
	\begin{theorem}
		The total number of positions of the form $"a, b, a + b"$ where all piles have $k$ or less digits (in binary form) is given by the formula $\frac{3^k}{2} - 2^{k}  + \frac{1}{2}$.
	\end{theorem}
	
	Proof. We can use the formula from Theorem \ref{sec:lemma-7.7} to calculate the total number of positions with two piles with $k$ digits. To compute the total number of positions where all piles have $k$ or less digits, the sum of a geometric series formula can be used. The exact working follows trivially.   
	
	\begin{definition}
		A \emph{transferable pair} is a pair of 0 and 1 chosen between the binary expansion of 2 piles. The digits must be chosen from the same digits place and the 1 must be from the smaller pile.
	\end{definition}
	
	\begin{lemma}
		The number of zero-nim positions that can be reached from a zero-nim position in 1 move is given by the total number of ways transferable pairs can be chosen between the first and second pile, second and third pile, and third and first pile.
	\end{lemma}	

	Proof. If the transfer of objects occur between transferable pairs, then the parity of 1s in each digit place of the binary expansion of the piles will remain the same. Thus, the nim sum will remain 0. 
	
	\begin{lemma}
		 The binary digits of the same digit place in all zero-nim positions are restricted to $(0, 1, 1), (1, 0, 1), (1, 1, 0), (0, 0, 0)$. Also, transfers between positions are restricted to moving 1s to 0s in the same digit place (ie.the binary digits of the same digit place is limited to the given 4 after transferring objects).
		 \label{sec:lemma-7.12}
	\end{lemma}

	Proof. For the first part of the lemma, we can see that these 4 sets of positions are the only positions with an even number of 1s.$\\$
	
	To prove the second part of the lemma, we see that if 1s are moved to 1s in the same digit place, then carrying of digits occur. With that in mind, there are 4 cases to consider for the next digits place after carrying of digits occur from $(1, 1, 0)$ (from second pile to first pile). $\\$
	
	Case 1. $(1, 1, 0) -> (0, 1, 0) \\$
	Case 2. $(1, 0, 1) -> (0, 0, 1) \\$
	Case 3. $(0, 1, 1) -> (1, 1, 1) \\$
	Case 4. $(0, 0, 0) -> (1, 0, 0) \\$	
	
	As the number of 1s in the next digits place is always odd, transfers between positions are restricted to moving 1s to 0s in the same digit place (note that moving 1s to 0s maintains the parity of 1s in that digits place).
	
	\begin{lemma}
		The total length of a transferable sequence of total pile sum $S$ is the total number of zero-nim positions with total pile sum $S$. The exact formula for the total number can be found with the formula given in Theorem \ref{sec:lemma-7.4}. 
	\end{lemma}

	Proof. 
	
	First, we note that all positions in a transferable sequence have the same sum and same nim-sum of 0.$\\$
	Also, by Lemma \ref{sec:lemma-7.12}, the binary digit places in all transferable positions are restricted to $(0, 1, 1), (1, 0, 1), (1, 1, 0), (0, 0, 0)$ and transfers between positions are restricted to moving 1s to 0s in the same digit place. $\\$
	
	Thus, all zero-nim positions with a sum $S$ are part of the same transferable sequence. This is because from any zero-nim position, any other zero-nim position of the same sum can be reached by reordering the binary digits of the same digit place. The only thing to note is that the direction of transfer between 2 transferable positions depends on the size of the piles before the transfer takes place.

	\section{Computer Science Considerations}
	
	In this section, we present computer programs for Theorems \ref{sec:lemma-7.2}, \ref{sec:lemma-7.3}, \ref{sec:lemma-7.4}, \ref{sec:lemma-7.5}. The programs can also be found in this Github repository: \url{https://github.com/ethank11k/Sharing-Nim} 
	$\\$
	
	\begin{lstlisting}[language=C++, caption={C++ code for Theorem \ref{sec:lemma-7.2} and Theorem \ref{sec:lemma-7.3}}]
#include <iostream>
using namespace std;
typedef long long ll;
ll sum, XOR_sum, x, bit, total_num = 1;
	
long long f() {
	if ((sum - XOR_sum) & 1) {
		cout << 0;
		return 0;
	}
	x = (sum - XOR_sum) / 2;
	for (ll i = 0; i < 50; i++) {
		bit = (1LL << i);
		if (x & bit) {
			if (XOR_sum & bit) {
				total_num = 0;
				break;
			}
		}
		else {
			if (XOR_sum & bit)
				total_num *= 2LL;
		}
	}
	if (total_num > 0 && x == 0) total_num -= 2;
	cout << total_num/2;
}
	
int main() {
	cin >> sum >> XOR_sum;
	f();
}
	\end{lstlisting}
	
	\begin{remark}
		The above has constant time complexity: $O(1)$.
	\end{remark}
	
	\begin{lstlisting}[language=C++, caption={C++ code for Theorem \ref{sec:lemma-7.4}}]
#include <iostream>
using namespace std;
typedef long long ll;
ll x, S, bit, temp_total, total_num = 0;

ll f(ll XOR_sum,ll sum) {
	temp_total = 1;
	if ((sum - XOR_sum) & 1) {
		return 0;
	}
	x = (sum - XOR_sum) / 2;
	for (ll i = 0; i < 50; i++) {
		bit = (1LL << i);
		if (x & bit) {
			if (XOR_sum & bit) {
				temp_total = 0;
				break;
			}
		}
		else {
			if (XOR_sum & bit)
				temp_total *= 2LL;
		}
	}
	if (temp_total > 0 && x == 0) temp_total -= 2;
	total_num += temp_total;
}


int main() {
	cin >> S;
	for(ll i=1;i<=S;i++) {
		f(i,S-i);
	}
	cout << total_num/6;
}

	\end{lstlisting}
	
	\begin{remark}
	The above has linear time complexity: $O(n)$.
	\end{remark}

	\begin{lstlisting}[language=C++, caption={C++ code for Theorem \ref{sec:lemma-7.5}}]

#include <iostream>
using namespace std;
typedef long long ll;
ll x, Q, S, bit, temp_total, total_num = 0;

ll f(ll XOR_sum,ll sum) {
	temp_total = 1;
	if ((sum - XOR_sum) & 1) {
		return 0;
	}
	x = (sum - XOR_sum) / 2;
	for (ll i = 0; i < 50; i++) {
		bit = (1LL << i);
		if (x & bit) {
			if (XOR_sum & bit) {
				temp_total = 0;
				break;
			}
		}
		else {
			if (XOR_sum & bit)
				temp_total *= 2LL;
		}
	}
	if (temp_total > 0 && x == 0) temp_total -= 2;
	total_num += temp_total;
}


int main() {
	cin >> Q;
	for(ll i=1;i<=Q;i++) {
		for(ll j=1;j<=i;j++) {
			f(j,i-j);
		}
	}
	cout << total_num/6;
}

	\end{lstlisting}
	
	\begin{remark}
		The above has quadratic time complexity: $O(n^2)$.
	\end{remark}
	
	\section{Nim inside a Nim}
	
	In this section, we will examine a specific case of Sharing Nim and identify both a winning strategy and winning/losing positions.
	
	The specific case to be examined are starting game positions of the form $1, a, b$, where $a$ and $b$ are any positive integer.
	
	\begin{lemma}
		If the absolute value of the difference between the second and third pile is greater than 1, the first pile with 1 object will not be touched by either player. 
		\label{sec:lemma-9.1}
	\end{lemma}

	Proof. If the absolute value of the difference between the second and third pile equals 1, the current player will choose to transfer the first pile to the lesser of the second or third pile. Then the game is equivalent to a game of Nim with 2 piles and the  player who made that move will win the game. $\\$
	
	However, if the absolute value of the difference between the second and third pile is greater 1, the current player will not touch the first pile. If the current player chooses to transfer or remove the first pile, then the game is once again equivalent to a game of Nim with 2 piles. Since the nim-sum would be greater than 0, the player who made that move will lose the game.
	
	\begin{lemma}
		The position $1, 2, 4$ is a $P$ position.
	\end{lemma}
	
	Proof. From this position, the next player's moves are limited to the following: $(a,b), (a, a, b), (a, b, a + b)$, where $a$ and $b$ are positive integers $\\$.
	
	All of these positions are clearly $N$ positions and thus the position $1, 2, 4$ is a $P$ position.
	
	\begin{theorem}
		Starting game positions of the form $1, a, b$, where the absolute value of the difference between the second and third pile is greater than 1, can be thought of as a game of 2-pile Nim.
		\label{sec:lemma-9.3}
	\end{theorem}

	Proof. First, by Lemma \ref{sec:lemma-9.1}, note that the first pile will not be touched by either player. Hence, only the second and third piles will be altered throughout the game.
	
	Since the position $1, 2, 4$ is a $P$ position, the game can then be thought of as a game of Nim with the piles $a-2, b-4$ and whichever player finishes the turn with both piles with 0 objects will win the game. 
	
	This is because by Lemma \ref{sec:lemma-5.2}, all moves a player can make in 2 pile sharing Nim maintains the property of Nim (all strategies of Nim apply as well).
	
	\begin{theorem}
		The winning strategy is to finish every move with the nim-sum of piles $a-2, b-4$ equal to 0.   
		\label{sec:lemma-9.4}
	\end{theorem}

	Proof. As shown in Lemma \ref{sec:lemma-5.2}, from a position where the nim-sum is not equivalent to 0, the next move can always be done so that the nim-sum is equivalent to 0. On the other hand, if the nim-sum is equivalent to 0, the next move will always result in a non-zero nim-sum.
	
	We now examine a generalized version of the above winning position in 3 pile Sharing Nim.
	
	\begin{lemma}
		A starting game position of the form $a, b, b + 2$, where $a$ and $b$ are positive integers, is a winning position. 
	\end{lemma}

	From this position, the first player can remove $a-1$ objects from the first pile, forming the position $1, b, b+2$. Then by Theorems \ref{sec:lemma-9.3} and \ref{sec:lemma-9.4}, the game is equivalent to a game of Nim. Since the Nim-sum of the 2 piles $b-2, b-2$ is 0. the first player will win the game.

	\section{Further Questions}

	Much remains to be learned about the game of Sharing Nim, especially in the general case. The author hopes that the study of Sharing Nim and its general winning strategy will continue. 
	
	\bibliographystyle{amsplain}

\end{document}